\documentclass[11pt,reqno]{amsart}
\usepackage{graphicx,amscd}
\usepackage{amsfonts,amsmath,amssymb}
\newcommand{\rl}{{\mathbb{R}}}
\newcommand{\cx}{{\mathbb{C}}}

\newcommand{\om}{\omega}

\newcommand{\dbar}{\overline{\partial}}

\newcommand{\e}{\varepsilon}

\newcommand{\dist}{{\mathrm{dist}}}

\newcommand{\B}{\mathbb{B}}

\newtheorem{theorem}{Theorem}
\newtheorem{lemma}[theorem]{Lemma}

\theoremstyle{definition}

\begin{document}
\title{Lagrangian inclusion with an open Whitney umbrella is rationally convex}
\author{Rasul Shafikov* and Alexandre Sukhov**}
\begin{abstract}
It is shown that a Lagrangian inclusion of a real surface in $\cx^2$ with a standard open Whitney umbrella
and double transverse self-intersections is rationally convex.
\end{abstract}

\maketitle

\let\thefootnote\relax\footnote{MSC: 32E20, 53D12. 
Key words: rational convexity, Lagrangian submanifolds
}

* Department of Mathematics, the University of Western Ontario, London, Ontario, N6A 5B7, Canada,
e-mail: shafikov@uwo.ca. The author is partially supported by the Natural Sciences and Engineering 
Research Council of Canada.

**Universit\'e des Sciences et Technologies de Lille, 
U.F.R. de Math\'ematiques, 59655 Villeneuve d'Ascq, Cedex, France,
e-mail: sukhov@math.univ-lille1.fr. The author is partially supported by Labex CEMPI.

%%%%%%%%%%%%%%%%%%%%%%%% section %%%%%%%%%%%%%%%%%%%%%%%
\section{Introduction}

This paper is concerned with the study of rational convexity of compact real surfaces
in $\cx^2$. A compact set $X$ in $\cx^n$ is {\it rationally convex} if for every point
$p$ in the complement of $X$ there exists a complex algebraic hypersurface passing
through $p$ and avoiding~$X$. See Stout~\cite{St} for a comprehensive treatment of this 
fundamental notion.

A nondegenerate closed 2-form $\omega$ on $\cx^2$ is called a {\it symplectic form}. 
By Darboux's theorem every symplectic form is locally equivalent to the standard form
$$
\omega_{\rm st} = \frac{i}{2} (dz \wedge d\bar z + dw \wedge d\bar w) = dd^c\, \phi_{\rm st}, \ \ 
\phi_{\rm st}= |z|^2 + |w|^2,
$$
where $(z,w)$, $z = x + iy$, $w = u + iv$ are complex coordinates in $\cx^2$, and 
$d^c =i(\overline\partial -\partial)$.  If a symplectic form $\om$ is of bidegree $(1,1)$ and strictly 
positive, it is called a {\it K\"ahler form}. 
A strictly plurisubharmonic function $\phi$ is called a potential of $\om$ if $dd^c \phi = \om$.
A real $n$-dimensional submanifold $S\subset \cx^n$ is called {\it Lagrangian} for $\om$ 
if $\omega|_S=0$. According to a  theorem of Duval and Sibony~\cite{DS}, a compact 
$n$-dimensional submanifold of $\cx^n$ is rationally convex if and only if it is Lagrangian for 
some K\"ahler form. This result displays a connection between rational convexity and 
symplectic properties of real submanifolds.

Being Lagrangian imposes certain topological restrictions on a submanifold, for example, 
the only compact orientable surface that admits a Lagrangian embedding into
$(\cx^2,\omega_{\rm st})$ is a torus. On the other hand, according to the result of Givental~\cite{Giv},
any compact surface (orientable or not) admits a {\it Lagrangian inclusion} into $\cx^2$, i.e., a 
smooth map $\iota: S \to \cx^2$ which is a local Lagrangian embedding except a finite set of singular points
that are either transverse  double self-intersections or the so-called {\it open Whitney umbrellas}. 
The {\it standard} open Whitney umbrella is a map 
\begin{eqnarray}\label{st-umb}
\pi: \rl^2_{(t,s)} \ni (t,s) \mapsto \left(ts,\frac{2t^3}{3},t^2,s\right) \in \rl^{4}_{(x,u,y,v)}.
\end{eqnarray}
The open Whitney umbrella is then defined as the image of the standard umbrella under a local symplectomorphism, 
i.e., a local diffeomorphism that preserves the form $\omega_{\rm st}$. It was proved by Gayet~\cite{G}
that an immersed Lagrangian (with respect to some K\"ahler form) submanifold in $\cx^n$ with 
transverse double self-intersections is also rationally convex. This was generalized to certain nontransverse 
self-intersections by Duval and Gayet~\cite{DG}.

The goal of this paper is show how the technique of~\cite{DS}, \cite{G}, and~\cite{DG} can be adapted to
prove rational convexity of a Lagrangian inclusion with one standard open Whitney umbrella. More precisely,
we prove the following.

\begin{theorem}\label{t.1}
Let $\iota : S \mapsto (\cx^2,\om_{\rm st})$ be a Lagrangian inclusion of a compact surface $S$. 
Suppose that the singularities of $\iota$ consist of transverse double self-intersections and
one standard open Whitney umbrella. Then $\iota(S)$ is rationally convex in $\mathbb C^2$.
\end{theorem}

We remark that the standard open Whitney umbrella can be replaced by its image under a complex affine map
that preserves the symplectic form $\om_{\rm st}$.  The existence of Lagrangian inclusions satisfying the 
conditions of Theorem~\ref{t.1} follows from a recent result of Nemirovski and  Siegel~\cite{NS}.

%%%%%%%%%%%%%%%%%%%%%%%%%%%%%%%%%%%%%%%%%%%%%%%%%%%%
\section{Proof of Theorem~\ref{t.1}.}

We will identify $S$ and $\iota(S)$ as a slight abuse of notation. The ball of radius $\e$ centred at a point $p$ 
is denoted by $\B(p,\e)$, and  the standard Euclidean distance between a point $p\in \cx^n$ and a 
set $Y\subset \cx^n$ is denoted by $\dist(p,Y)$. Our approach is a modification of the method of Duval-Sibony 
and Gayet. The main tool here is the following result. 

\begin{lemma}[\cite{DS}, \cite{G}] \label{l.g}
Let $\phi$ be a plurisubharmonic $C^\infty$-smooth function on $\cx^n$, and let $h$ be a $C^\infty$-smooth function
on $\cx^n$ such that 

\begin{itemize}
	\item[(1)] $|h| \le e^\phi$, and $X:= \{|h|=e^\phi \}$ is compact;
	\item[(2)] $\dbar h = O (\dist(\cdot, S)^{\frac{3n+5}{2}})$;
	\item[(3)] $|h| = e^\phi$ with order 1 on $S$; 
	\item[(4)] For any point $p\in X$ at least one of the following conditions hold: (i) $h$ is holomorphic in a neighbourhood 
	of $p$, or (ii) $p$ is a smooth point of $S$, and $\phi$ is strictly plurisubharmonic at $p$.
\end{itemize}

Then $X$ is rationally convex. 
\end{lemma}

The proof of Theorem~\ref{t.1} consists of finding the functions $\phi$ and $h$ that satisfy Lemma~\ref{l.g}
and such that the set $X$ contains $S$ and is contained in the union of $S$ with the balls of arbitrarily small radius
centred at singular points of $S$. This will be 
achieved in three steps: first we construct a closed $(1,1)$-form $\om$ that vanishes near singular points of $S$ and such 
that $\om|_S=0$. This is done is Section~\ref{s.om}. The form $\om$ is a modification of the standard symplectic form 
$\om_{\rm st}$ in $\cx^2$ near singular points of $S$. Near self-intersection points this is done in the paper of 
Gayet~\cite{G}, and so we will deal with the umbrella point. Secondly, from $\om$ and its potential $\phi$ we construct 
the required function $h$. This is done in Section~\ref{s.om2}. 
In the last step, in Section~\ref{s.om3}, we replace $\phi$ with a function $\phi + \rho$, for a 
suitable $\rho$, so that the pair $\{\phi+\rho, h\}$ satisfies all the conditions of Lemma~\ref{l.g}.

%%%%%%%%%%%%%%%%%%%%%%%%%%%%%%%%%%%%%%%%%%%
\subsection{The form $\om$.}\label{s.om}
Near the umbrella point the Lagrangian inclusion map $\iota$ coincides with $\pi$ given by~\eqref{st-umb}. 
For a function $f$ we have 
$$
d^c f = -f_y dx + f_x dy - f_v du + f_u dv.
$$ 
Direct computations show that
$\pi^*d^c \phi_{\rm st} = -2t^2 s dt - \frac{2}{3}t^3 ds$. Consider the pluriharmonic function
$\zeta = \frac{v^2}{2} - \frac{u^2}{2}$.
Then $\pi^*d^c \zeta = \pi^*d^c \phi_{\rm st}$. The function
$$
\phi = \phi_{\rm st} - \zeta 
$$ 
is strictly plurisubharmonic and satisfies 
\begin{equation}\label{e.*}
\pi^*d^c \phi=0.
\end{equation}
Let $r: \rl^{+} \to \rl^+$ be a smooth increasing convex function such that $r(t)=0$
when $t\le\e_1$ and $r(t) = t -c$ when $t>\e_2$, for some suitably chosen $c>0$ and $0<\e_1<\e_2$. 
We choose $\e_2>0$ so small that the set $\{\phi<\e_2\}$ does not contain singular points
of $S$ except the origin. Let
\begin{equation}\label{e.ome}
\om = dd^c(r \circ \phi).
\end{equation}
Then
$\pi^*\omega =0$
by~\eqref{e.*}. Therefore, the surface $S$ remains Lagrangian with respect to the form $\omega$. This gives us 
the required modification of $\om_{\rm st}$. By construction there exist two neighbourhoods $U\Subset U'$ of the
origin such that $\om|_U=0$ and $\omega =\omega_{\rm st}$ in $\cx^2 \setminus U'$, while the potential
changed globally.

Denote by $p_1, \dots, p_N$ the points of self-intersection of $S$, and by $p_0$ the standard umbrella. 
Then \cite[Prop. 1]{G} gives further modification $\tilde \omega$ of the form $\omega$  in~\eqref{e.ome}, 
near the self-intersection points. Combining everything together yields the following result.

\begin{lemma}\label{l.omega}
Given $\e>0$ sufficiently small, there exists a $(1,1)$-form $\tilde \om$ and $\e_1>0$, such that 
\begin{itemize}
	\item[(i)] $\tilde \om|_S = 0$;
	\item[(ii)] $\tilde \om = \om$ on $\cx^2 \setminus \cup_{j=0}^N \B(p_j, \e)$.
	\item[(iii)] $\tilde \om$ vanishes on  $\mathbb B(p_j,\e_1)$, $j=0,\dots,N$.
	
\end{itemize}
Furthermore, there exists a smooth function $\tilde \phi$ on $\cx^2$ such that $dd^c \tilde \phi = \tilde \omega$.
The function $\tilde \phi$ is plurisubharmonic on $\mathbb C^2$, and strictly plurisubharmonic on 
$\cx^2 \setminus \cup_{j=0}^N \B(p_j, \e)$.
\end{lemma}

%%%%%%%%%%%%%%%%%%%%%%%%%%%%%%%%%%%%%%
\subsection{The function $h$.}\label{s.om2}
Let $\iota: S \to \cx^2$ be a Lagrangian inclusion, and $\tilde\phi$ be the potential of the form $\tilde\omega$ 
given by Lemma~\ref{l.omega}. For simplicity we drop tilde from the notation.
In this subsection we recall the construction in \cite{DS} and \cite{G} of a smooth function $h$ on $\cx^2$ such that 
$|h|\left|_S = e^\phi \right.$ and $\dbar h(z) = O(\dist(z,S)^6)$.  The two conditions, that 
$\dbar h$ vanishes on $S$ and that $\phi - \log |h|$ vanishes on $S$ with order 1 imply that 
$\iota^*(d^c \phi - d({\rm arg}h))=0$. The latter condition 
can be met by further perturbation of  $\phi$. 

Let $\tilde S$ be the deformation retract of $S$. Note that it exists because near the umbrella 
point the surface $S$ is the graph of a continuous vector-function. Let $\gamma_k$, $k=1,\dots,m$, be the basis in
$H_1(\tilde S, \mathbb Z) \cong H_1(S, \mathbb Z)$ supported on $S$. Using de Rham's theorem 
and an argument similar to that of Lemma~\ref{l.omega} one can find smooth functions $\psi_k$ 
with compact support in $\cx^2$ that vanish on $S\cup(\cup_j B(p_j,\e))$, where $B(p_j,\e)$ 
are the balls around the singular points on $S$ as in Lemma~\ref{l.omega}, such that 
$\int_{\gamma_k}\iota^* d^c\psi_l = \delta_{kl}$.
Further, one can find small rational numbers $\lambda_k$ and an integer $M$, such that for the function 
\begin{equation}\label{e.phi-til}
\tilde \phi = M \left(\phi+ \sum_{j=1}^m \lambda_k \psi_k \right)
\end{equation}
the form $\iota^* d^c\tilde\phi$ is closed on $S$ and has periods which are multiples of $2\pi$. Then
there exists a $C^\infty$-smooth function $\mu: S \to \rl/2\pi\mathbb Z$ that vanishes on the intersection
of $S$ with $B(p_j,\e)$, $j=0,\dots,N$, and such that $\iota^* d^c \tilde \phi = d\mu$. By~\cite{HW}, there exists a 
function $h$ defined on $\cx^2$ such that 
$$
h|_S = e^{\tilde \phi + i\mu}|_S
$$
and $\dbar h(z) = O(\dist(z,S)^6)$. It follows that $\tilde \phi - \log|h|$ vanishes to order $1$ on $S$.
Note that $h$ is constant near singular points of $S$. 

%%%%%%%%%%%%%%%%%%%%%%%%%%%%%%%%%
\subsection{The function $\phi$.}\label{s.om3}
Again, for simplicity of notation we denote by $\phi$ the function~\eqref{e.phi-til} constructed in 
Section~\ref{s.om2}. It does not yet satisfy the conditions of Lemma~\ref{l.g}
because there are still some smooth points on $S$ where the function $h$ is not holomorphic
and $\phi$ is not strictly plurisubharmonic. For this we will replace $\phi$ by a 
function $\tilde \phi = \phi + c\cdot\rho$, where the function $\rho$ will be constructed using local 
polynomial convexity of $S$, and $c>0$ will be a suitable constant.

We recall our result from~\cite{SS1,SS2}.

\begin{lemma}\label{l.FS} 
Let $S$ be a Lagrangian inclusion in $\cx^2$, and let $p_0,\dots, p_N$ be its singular points. 
Suppose that $S$ is locally polynomially convex near every singular point. Then there exists a neighbourhood 
$\Omega$ of $S$ in $\cx^2$ and a continuous non-negative plurisubharmonic function $\rho$  on $\Omega$ 
such that $S\cap \Omega = \{ p \in \Omega: \rho(p) = 0 \}$. Furthermore, for every $\delta > 0$ one can 
choose $\rho = (\dist(z,S))^2$ on $\Omega \setminus \cup_{j=0}^N \B(p_j,\delta)$; in particular, it is 
smooth and strictly plurisubharmonic there.
\end{lemma}

The standard open Whitney umbrella is  locally polynomially convex by \cite{SS1}, and $S$ is locally polynomially
convex near transverse double self-intersection points by~\cite{SS2}. For the proof of the lemma
we refer the reader to~\cite{SS2}.

To complete the construction of the function $\phi$, we choose the function $\rho$ in Lemma~\ref{l.FS} 
with $\delta>0$ so small that the balls $\mathbb B(p_j, \delta)$ are contained in balls $\mathbb B(p_j, \e_1/2)$ 
given by Lemma~\ref{l.omega}. Note that $\rho$ is defined only in a 
neighbourhood $\Omega$ of $S$, but we can extend it as a smooth function with compact support in $\cx^2$. 
Consider now the function 
$$
\tilde \phi = \phi + c\cdot \rho .
$$
We choose the constant $c>0$ so small that the
function $\tilde \phi$ remains to be plurisubharmonic on $\cx^2$. At the same time, since $c>0$ and $\rho$ is
strictly plurisubharmonic on $S$ outside small neighbourhoods of singular points, we conclude that the function 
$\tilde \phi$ is strictly plurisubharmonic outside the balls $\mathbb B(p_j, \delta)$. 

The pair $\tilde\phi$ and $h$ now satisfies all the conditions of Lemma~\ref{l.g}. This completes the proof
of Theorem~\ref{t.1}.

%%%%%%%%%%%%%%% Biblio %%%%%%%%%%%%%%%%%%%%%%%

\end{document}